# Differential Game Model of Dispersed Material Drying

Malafeyev Oleg[1, a)] and Rylow Denis[1, b)]
and Zaitseva Irina[2, 3, c)] and Novozhilova Lydia[1, d)]

[1]*St. Petersburg State University, Faculty of Applied Mathematics and Control Processes,
7/9 Universitetskaya nab., St. Petersburg, 199034, Russia*
[2]*Stavropol State Agrarian University, Zootekhnicheskiy lane 12, Stavropol, 355017, Russia*
[3]*Stavropol branch of the Moscow Pedagogical State University, Dovatortsev str. 66 g, Stavropol, 355042, Russia*
[3]*St. Petersburg State University, Institute of Earth Sciences, Departments of Ecological Geology,
7/9 Universitetskaya nab., St. Petersburg, 199034, Russia*

[a)]malafeyevoa@mail.ru
[b)]Corresponding author: denisrylow@gmail.com
[c)]ziki@mail.ru
[d)]xopctt@rambler.ru

**Abstract.** Continuous and discrete game-theoretic models of dispersed material drying process are formalized and studied in the paper. The existence of optimal drying strategies is shown through application of results from the theory of differential games and dynamic programming. These optimal strategies can be found numerically.

## INTRODUCTION

Physical processes involving energy or material transfer, including mechanical movement, have long been described with the help of evolution equations, especially differential ones. This is the standard method of modelling electromagnetic systems, including electromechanical and computer components, power units, grids, antennas, etc. [1-6]. Such framework also allows to set and research the problems of control, stability and reliability [7-10].

Another less popular approach to modelling of physical processes is to use the game theory, which is mainly used in studies of socioeconomic systems and cooperation [11-23]. Examples of game theory application to the problems of physics can be found in [24-33].

In this paper, the dispersed material drying process is formalized as a two player game between the "nature" and a drying apparatus operator. The paper is organized as follows. In section "the continuous model of dispersed material drying process" we introduce the general framework of the model and study the optimal strategies of both players. In section "the discrete model of dispersed material drying process" the framework introduced in the previous section is extended to the case of $n$-step drying process. The conclusions are drawn in the last section.

## DIFFERENTIAL GAME MODEL OF DISPERSED MATERIAL DRYING PROCESS

Measurable function $t(\tau), t:[0, J^*] \to T$ is an admissible control of the operator. Here $[0, J^*]$ is the time interval on which the drying process takes place and $T \subset [0, \infty)$ is a compact set of all possible temperatures the dispersed materials can be dried at in the apparatus. The measurable function $\alpha(\tau), \alpha:[0, J^*] \to Q$ is an admissible control of the "nature", where $Q$ is a compact set of all possible parameters relevant to the drying process that can not be controlled by the operator. Denote by $x$ the humidity of the dispersed material. Trajectory of the process, corresponding to the pair of controls $(t,\alpha)$, is an absolutely continuous function $x:[0,J^*] \to R_1^+$ such that for almost all $\tau:[0,J^*]$ the equation $\dot{x}(\tau) = f(x(\tau),t(\tau),\alpha(\tau))$ holds.

It is assumed that the drying process proceeds over the interval $[0,J^*]$ in such a way that the trajectory reaches the given set $\tilde{X} \subset R_1^+$ called the terminal set of the game before or at the time moment $J^*$.

Let the energy cost of the drying process be characterized by a given function $e(t,\alpha,x,\tau)$. If $t(\tau)$, $\alpha(\tau)$ are control functions of the operator and the "nature" over the interval $[0,J^*]$, $x(\tau)$ is the trajectory of the drying process corresponding to them, then $\int_0^\tau e(t(\tau),\alpha(\tau),\tau)d\tau = E(t,\alpha,x(t,\alpha))$ is energy cost for the drying process.

Within the framework of our model the following two problems are formalized:
1. to find the control function $t(\tau)$ – thermal regime of the drying process – that brings the humidity of the dispersed material from the initial conditions to the final conditions $x(\tau) \in \tilde{X}$ with minimal energy cost under all possible external parameters – controls of the "nature";
2. to find the control function $t(\tau)$ that brings the humidity of the dispersed material to the final conditions $x(\tau) \in \tilde{X}$ in the minimal time under all possible external parameters.

For definiteness the first problem is considered. It is natural to consider the upper (majorant) game $\Gamma$ (cf. [23, 29]) within the framework of the problem. Then the piecewise-program strategies are defined as follows. The strategy of the minimizing operator $\overline{T}$ is a pair $\overline{T} = (\sigma',\varphi_{\sigma'})$, where $\sigma'$ is the finite partition of the interval $[0,J^*]$, $\varphi_{\sigma'}$ is the lower strategy corresponding to the partition $\sigma'$, that is the mapping which associates the informational status of the operator at the time moment $\tau_l \in \sigma'$ with the measurable control $t_l(\cdot)$ over the interval $[\tau_l,\tau_{l+1}]$.

The strategy of the "nature" $\overline{\alpha}$ is the set $\{\psi_\sigma\}_{\sigma \in \Sigma}$, where $\Sigma$ is the set of all finite partitions of the interval $[0,J^*]$, $\psi_\sigma$ is the upper strategy of the "nature" in the multi-step upper game $\Gamma^{\sigma'}$ corresponding to the partition $\sigma'$, that is the mapping which associates the informational status of the "nature" at the time moment $\tau_l \in \sigma'$ with the measurable control $\alpha_l(\cdot)$ on the interval $[\tau_l,\tau_{l+1}]$. Recall that in upper game $\overline{\Gamma}^\sigma$ the operator knows the state of the process at any time moment $t_k \in \sigma$. The "nature" knows the state of the process at any time moment $t_k \in \sigma$ as well, along with the control of the operator at the interval $[t_k,t_{k+1}]$. Having found the pair of strategies $(\overline{T},\overline{\alpha})$, one can construct the single set of controls $t(\cdot),\alpha(\cdot)$ and the corresponding trajectory of the process $x(\tau) = \chi(\overline{T},\overline{\alpha})$ over the whole interval $[0,J^*]$. It is known [23], that there exists $\varepsilon$-saddle point for all $\varepsilon > 0$ in the game $\overline{\Gamma}$, that is there exists such a pair of strategies $\overline{T}^\varepsilon,\overline{\alpha}^\varepsilon$, that for all $\overline{T} \in T, \overline{\alpha} \in Q$ the inequality $E\left(\overline{T}^\varepsilon,\overline{\alpha},\chi\left(\overline{T}^\varepsilon,\overline{\alpha}\right)\right) + \varepsilon \le E\left(\overline{T}^\varepsilon,\overline{\alpha}^\varepsilon,\chi\left(\overline{T}^\varepsilon,\overline{\alpha}^\varepsilon\right)\right) \le E\left(\overline{T},\overline{\alpha}^\varepsilon,\chi\left(\overline{T},\overline{\alpha}^\varepsilon\right)\right) - \varepsilon$ holds. The value $E\left(\overline{T}^\varepsilon,\overline{\alpha}^\varepsilon,\chi\left(\overline{T}^\varepsilon,\overline{\alpha}^\varepsilon\right)\right) - \varepsilon$ is the upper limit of the energy needed to dry the dispersed material regardless of the strategy the "nature" plays. Here $E(\cdot)$ is payoff in the situation $(\overline{T}^\varepsilon,\overline{\alpha}^\varepsilon)$.

## THE DISCRETE MODEL OF DISPERSED MATERIAL DRYING PROCESS

Assume that the dispersed material is to be processed in a fixed sequence $i_n$ in drying apparatuses $i = 1,2,\ldots,n$. At the $i$-th step $(i = \overline{1,n})$ the drying process consumes the amount of energy $e_i(t_i,\alpha_i,x_{i-1})$ which depends on the thermal regime chosen by the operator and the vector of external parameters chosen by the "nature". The temperature is assumed be constant at each $i$-th step. During the whole drying process the amount of energy

$$\sum_{i=1}^{n} e_i(t_i, \alpha_i, x_{i-1}) = \varepsilon_n(t_n, \alpha_n) = \varepsilon_n(t_1,...,t_n, \alpha_1,...,\alpha_n) = \varepsilon_n(t^n, \alpha^n)$$ is consumed in $n$ steps. At the each $i$-th step the limits of possible temperature variations and the vectors of external parameters are given: $t_i \in T_i = \left[t_i^1, t_i^2\right]; \alpha_i \in Q_i = \left[\alpha_i^1, \alpha_i^2\right]$. As in the previous section, it is necessary to find such sequence of controls $t = \{t_i\}$, which would guarantee the lowest energy cost under all possible external conditions $\alpha = \{\alpha_i\}$. Note that the material can be placed in the $i$-th drying apparatus long enough for external parameters to change. Therefore in this case, the set of external parameters vectors is a set of functions on some time interval.

This set of functions is taken to be a compact space, denote it by $Q$. Assume that the number of steps $n$ is such that the terminal set $\tilde{X}$ is reached for any control $\{t_i\}$ and any functions of external parameters $\{\alpha_i\}$.

Let us consider the $n$-step drying process for any external parameters with the initial humidity is $x_0$: $F_n(x_0) = \min_t \max_\alpha e_n(t_1^n, \alpha^n, x^n)$. Here $F_n(x_0)$ is the minimum amount of guaranteed energy consumed drying the drying process.

It is easy to show that principles of dynamic programming can be applied here: it is possible to construct recurrent equations, where two consecutive values of the Belmann function are connected: $F_n(x_0) = \min_{t_1} \max_{\alpha_1} \{e_1(t_1, \alpha_1, x_0) - F_{n-1}(x_1)\}$, $x_1 = x_0 + \Delta t_1 f(x_0, t_1, \alpha_1)$.

## CONCLUSIONS

The dispersed material drying process is formalized as a continuous two-player game-theoretic model, the discrete model is considered as well. An algorithm for finding optimal strategies in the model is given. The results obtained in this paper can be used in further studies of the drying process including numerical computation of optimal strategies.